\documentclass[11pt]{article}

\usepackage{latexsym}
\usepackage{amssymb}

\newtheorem{theorem}{Theorem}
\newtheorem{corollary}{Corollary}

\newtheorem{remark}{Remark}

\begin{document}
{
\begin{center}
{\Large\bf
On extensions of commuting tuples of symmetric and isometric operators.}
\end{center}
\begin{center}
{\bf S.M. Zagorodnyuk}
\end{center}

\section{Introduction.}

Extension problems for commuting tuples of symmetric and isometric operators have a long history, see, 
e.g.,~\cite[Chapter VIII]{cit_7000_Berezanskii}, \cite{cit_8900_Coddington},
\cite{cit_580_Phillips},
\cite{cit_8300_Ismagilov},
\cite{cit_8350_Slinker},
\cite{cit_8500_Jorgensen}, 
\cite{cit_590_Kochubei},
\cite{cit_8700_Schmuedgen},
\cite{cit_8800_Bohonov},
\cite{cit_8800__Vasilescu},
\cite{cit_8800_Albrecht_Vasilescu},
\cite{cit_8900__Vasilescu},
\cite{cit_10000_Zagorodnyuk} and references therein.
A set
\begin{equation}
\label{f1_10}
\mathcal{T} = (A_1,...,A_\rho, B_1,...,B_\tau),\qquad \rho,\tau\in\mathbb{Z}_+,\ \rho+\tau\geq 2,
\end{equation}
where all $A_j$ are symmetric operators and all $B_k$ are linear isometric operators,
having a joint invariant dense domain $\mathcal{D}$ in a Hilbert space $H$ and pairwise commuting,
is said to be \textit{a commuting tuple of symmetric and isometric operators (with a joint invariant dense domain) in $H$}.
The case $\rho=0$ (or $\tau=0$) means that operators $A_j$ (respectively $B_k$) are absent.  

Commuting tuples of operators appear in a natural way when dealing with various moment problems, see, e.g.~\cite{cit_7000_Berezanskii},
\cite{cit_2000_Fuglede}.
In this context there appears a natural question:

\noindent
\textbf{Question 1}: \textit{Do there exist self-adjoint operators $\widetilde A_j\supseteq A_j$ and
unitary operators $\widetilde B_k\supseteq B_k$ in a Hilbert space $\widetilde H\supseteq H$, all pairwise commuting?}

In~\cite{cit_9000_Zagorodnyuk_MFAT}, for the case $\rho=2$, $\tau=0$, we present some conditions which ensure the existence of
such commuting (self-adjoint) extensions in the same space $H$. Namely, the following theorem was obtained:

\begin{theorem} (\cite[Theorem 1]{cit_9000_Zagorodnyuk_MFAT})
\label{t1_1}
Let $\mathbf{A}$ be a symmetric operator and $\mathbf{B}$ be an essentially self-adjoint operator
with a common domain $\mathcal{D} = D(\mathbf{A}) = D(\mathbf{B})$
in a Hilbert space $\mathbf{H}$, $\overline{\mathcal{D}} = \mathbf{H}$, and
$$
\mathbf{A} \mathcal{D} \subseteq \mathcal{D},\quad \mathbf{B} \mathcal{D} \subseteq \mathcal{D}; $$
$$
\mathbf{A} \mathbf{B} = \mathbf{B} \mathbf{A}. $$
Suppose also that for some $z_0\in \mathbb{C}\backslash \mathbb{R}$, the operator $\mathbf{B}$
restricted to the domain $(\mathbf{A}-z_0 E_\mathbf{H})\mathcal{D}$ is essentially self-adjoint in a
Hilbert space
$(\overline{\mathbf{A}} - z_0 E_\mathbf{H}) D(\overline{\mathbf{A}})$.

If there exists a conjugation $\mathbf{J}$ in $\mathbf{H}$ such that
$\mathbf{J} \mathcal{D}\subseteq \mathcal{D}$, and
$$
\mathbf{A} \mathbf{J} = \mathbf{J} \mathbf{A},\quad \mathbf{B} \mathbf{J} = \mathbf{J} \mathbf{B}, $$
then there exists a self-adjoint operator $\widetilde{\mathbf{A}}\supseteq \mathbf{A}$, whichcommutes with $\overline{\mathbf{B}}$.
\end{theorem}

Recall that a conjugation $J$ in a Hilbert space $H$ is an antilinear operator, defined on the whole space $H$, such that:
$J^2 = E_H$, $(Jf,Jg)_H = (g,f)_H$, for all $f,g\in H$.
The proof of Theorem~\ref{t1_1} essentially used the remarkable Godi\v{c}-Lucenko Theorem (see~\cite{cit_2000_GL}, \cite{cit_3000_GP}).

Let us briefly describe the content of the present paper.
In this paper we aim to generalize our extension result to the case of commuting tuples of symmetric and isometric operators.
For this purpose we shall obtain a multidimensional extension of the Godi\v{c}-Lucenko Theorem (Theorem~\ref{t2_1}). 
The idea of its proof is similar to the proof of the Godi\v{c}-Lucenko Theorem in the paper of Garcia and Putinar~\cite{cit_3000_GP}.
It is crucial for our application that one of conjugations can be chosen to be the same for all unitary operators in a tuple.
Another crucial point is that the construction of an extension in the proof of Theorem~\ref{t1_1} used only one of two conjugations
provided by the Godi\v{c}-Lucenko Theorem.
In a consequence, we can generalize our result to several commuting symmetric operators.  Under some suitable conditions we add some
isometric operators to the commuting tuple as well~(Theorem~\ref{t2_3_m}). In the case of bounded $A_2,...,A_\rho$ the situation
simplifies considerably~(Corollary~\ref{c2_1})
We shall illustrate our results on a multidimensional moment problem~(Theorem~\ref{t2_4}).
The latter problem generalizes Devinatz's (or power-trigonometric) moment problem (see~\cite{cit_8950_Zagorodnyuk_JOT}).

\noindent
{\bf Notations. } We denote by $\mathbb{R},\mathbb{C},\mathbb{N},\mathbb{Z},\mathbb{Z}_+$
the sets of real numbers, complex numbers, positive integers, integers and non-negative integers,
respectively. By $\mathbb{Z}_{k,l}$ we mean all integers $j$ satisfying the following inequality:
$k\leq j\leq l$; ($k,l\in\mathbb{Z}$).
The Cartesian product of $n$ copies of a number set $S$ is denoted by $S^n$; $n\in\mathbb{N}$ (e.g. $\mathbb{Z}^n$, $\mathbb{Z}_+^n$, etc.).

Everywhere in this paper, all Hilbert spaces are assumed to be separable, and isometric operators are supposed to be linear. By
$(\cdot,\cdot)_H$ and $\| \cdot \|_H$ we denote the scalar product and the norm in a Hilbert space $H$,
respectively. The indices may be omitted in obvious cases.
For a set $M$ in $H$, by $\overline{M}$ we mean the closure of $M$ in the norm $\| \cdot \|_H$.
For
$\{ x_k \}_{k\in S}$, $x_k\in H$, by
$\mathop{\rm Lin}\nolimits \{ x_k \}_{k\in S}$ we mean the set of linear combinations of vectors $\{ x_k \}_{k\in S}$
and $\mathop{\rm span}\nolimits \{ x_k \}_{k\in S} =
\overline{ \mathop{\rm Lin}\nolimits \{ x_k \}_{k\in S} }$.
Here $S$ is an arbitrary set of indices.

The identity operator in $H$ is denoted by $E_H$. For an arbitrary linear operator $A$ in $H$,
the operators $A^*$,$\overline{A}$,$A^{-1}$ mean its adjoint operator, its closure and its inverse
(if they exist). By $D(A)$ and $R(A)$ we mean the domain and the range of the operator $A$.

Denote $D_{r,l} = \mathbb{R}^{r} \times [-\pi,\pi)^{l}
= \{ (x_1,x_2,...,x_r,\varphi_1,\varphi_2,...,\varphi_l),$
$x_j\in \mathbb{R},\ \varphi_k\in [-\pi,\pi),\ 1\leq j\leq r,\ 1\leq k\leq l \}$, $r,l\in \mathbb{Z}_+$.
For elements $u\in D_{r,l}$ we briefly write: $u = (x,\varphi)$, $x=(x_1,x_2,...,x_r)$,
$\varphi = (\varphi_1,\varphi_2,...,\varphi_l)$. We mean that $D_{r,0} = \mathbb{R}^{r}$;
$D_{0,l} = [-\pi,\pi)^{l}$.
By $\mathfrak{B}(D_{r,l})$ we mean the set of all Borel subsets of $D_{r,l}$.

\noindent
Let $M(\delta) = (m_{i,j}(\delta))_{i,j=0}^{N-1}$ be a
$\mathbb{C}_{N\times N}^\geq$-valued measure on $\mathfrak{B}(D_{r,l})$, and
$\tau = $ $\tau_M (\delta)$ $:=$ $\sum_{k=0}^{N-1} m_{k,k} (\delta)$;
$M'_\tau = $ $(m'_{k,l})_{k,l=0}^{N-1} = ( dm_{k,l}/ d\tau_M )_{k,l=0}^{N-1}$; $N\in \mathbb{N}$.
We denote by $L^2(M)$ a set of all (classes of the equivalence of)
vector-valued functions
$f: D_{r,l}\rightarrow \mathbb{C}_N$, $f = (f_0,f_1,\ldots,f_{N-1})$, such that
$$ \| f \|^2_{L^2(M)} := \int_{D_{r,l}}  f(u) \Psi(u) f^*(u) d\tau_M  < \infty. $$
As it is well known, the set $L^2(M)$ is a Hilbert space with the following scalar product:
$$ ( f,g )_{L^2(M)} := \int_{D_{r,l}}  f(u) \Psi(u) g^*(u) d\tau_M,\qquad f,g\in L^2(M). $$
Set
$$ W_n f(x,\varphi) = e^{i\varphi_n} f(x,\varphi),\qquad f\in L^2(M);\ 1\leq n\leq l; $$
and
$$ X_m f(x,\varphi) = x_m f(x,\varphi), $$
$$ f(x,\varphi)\in L^2(M):\ x_m f(x,\varphi)\in L^2(M);\
1\leq m\leq r. $$
Operators $W_n$ are unitary, while operators $X_m$ are self-adjoint.

\section{A multidimensional version of the Godi\v{c}-Lucenko Theorem and extensions of tuples of operators.}

We shall use a model for a set of commuting self-adjoint and unitary operators with the spectrum of 
a finite multiplicity provided by~\cite{cit_9500_Zagorodnyuk_M_Notes}. 
In order to state the corresponding result, we need to recall some relevant definitions from~\cite{cit_9500_Zagorodnyuk_M_Notes}.
Consider a set
\begin{equation}
\label{f2_1}
\mathcal{A} = (S_1,S_2,...,S_\mathbf{r},U_1,U_2,...,U_\mathbf{l}),\quad \mathbf{r},\mathbf{l}\in \mathbb{Z}_+:\
\mathbf{r}+\mathbf{l}\not=0,
\end{equation}
where $S_j$ are self-adjoint operators, $U_k$ are unitary operators in a separable Hilbert space $H$,
$1\leq j\leq \mathbf{r}$, $1\leq k\leq \mathbf{l}$.
In the case $\mathbf{r}=0$ operators $S_j$ are absent, while the case $\mathbf{l}=0$ means that $U_k$ disappear. 
The set $\mathcal{A}$ is said to be {\bf a $SU$-set of order $(\mathbf{r},\mathbf{l})$}.
The set $\mathcal{A}$ is called {\bf commuting}, if operators $S_j$,$U_k$ pairwise commute.
In the latter case there exists the spectral measure $E(\delta)$, $\delta\in \mathfrak{B}( D_{\mathbf{r},\mathbf{l}} )$, of $\mathcal{A}$ such that:

$$ S_j = \int_{ D_{\mathbf{r},\mathbf{l}} } x_j dE,\quad 1\leq j\leq \mathbf{r};\quad
U_k = \int_{ D_{\mathbf{r},\mathbf{l}} } e^{i\varphi_k} dE,\quad 1\leq k\leq \mathbf{l}. $$

A commuting $SU$-set $\mathcal{A}$ of order $(\mathbf{r},\mathbf{l})$
{\bf has the spectrum of order $d$}, if

\begin{itemize}
\item[1)] there exist vectors $h_0,h_1,...,h_{d-1}$ in $H$ such that
\begin{equation}
\label{f2_4}
h_i \in D(S_1^{m_1}S_2^{m_2}...S_\mathbf{r}^{m_\mathbf{r}}),\quad m_1,m_2,...,m_\mathbf{r}\in \mathbb{Z}_+,\
0\leq i\leq d-1;
\end{equation}
$$ \mathop{\rm span}\nolimits \{
U_1^{n_1} U_2^{n_2}...U_l^{n_\mathbf{l}} S_1^{m_1}S_2^{m_2}...S_\mathbf{r}^{m_\mathbf{r}} h_i, $$
\begin{equation}
\label{f2_5}
m_1,m_2,...,m_\mathbf{r}\in \mathbb{Z}_+;\ n_1,n_2,...,n_\mathbf{r}\in \mathbb{Z};\
0\leq i\leq d-1 \} = H;
\end{equation}

\item[2)] ({\it the minimality})
For an arbitrary $\widetilde d\in \mathbb{Z}_+:\ \widetilde d < d$, and an arbitrary
$\widetilde h_0,\widetilde h_1,...,\widetilde h_{\widetilde d-1}$ in $H$,
at least one of the conditions~(\ref{f2_4}),(\ref{f2_5}), with $\widetilde d$ instead of $d$, and
$\widetilde h_i$ instead of $h_i$, is not valid.

\end{itemize}

In the case $\mathbf{r}=0$, condition~(\ref{f2_4}) is redundant. Condition~(\ref{f2_5}) in the cases
$\mathbf{r}=0$, $\mathbf{l}=0$, does not have $U_k$ or $S_j$, respectively.
Set
$$ \vec e_i = (\delta_{0,i},\delta_{1,i},...,\delta_{N-1,i}),\qquad 0\leq i\leq N-1. $$

\begin{theorem}(\cite[Theorem 1]{cit_9500_Zagorodnyuk_M_Notes})
\label{t2_1_MN}
Let $\mathcal{A}$ be a commuting $SU$-set of order $(\mathbf{r},\mathbf{l})$ in a separable Hilbert space $H$,
having the spectrum of multiplicity $d$.
Let $x_0,x_1,...,x_{N-1}$, $N\geq d$, be elements of $H$ such that
\begin{equation}
\label{f2_9}
x_i \in D(S_1^{m_1}S_2^{m_2}...S_\mathbf{r}^{m_\mathbf{r}}),\quad m_1,m_2,...,m_\mathbf{r}\in \mathbb{Z}_+,\
0\leq i\leq N-1;
\end{equation}
$$ \mathop{\rm span}\nolimits \{
U_1^{n_1} U_2^{n_2}...U_l^{n_\mathbf{l}} S_1^{m_1}S_2^{m_2}...S_\mathbf{r}^{m_\mathbf{r}} x_i, $$
\begin{equation}
\label{f2_10}
m_1,m_2,...,m_\mathbf{r}\in \mathbb{Z}_+;\ n_1,n_2,...,n_\mathbf{r}\in \mathbb{Z};\
0\leq i\leq N-1 \} = H.
\end{equation}
Set
\begin{equation}
\label{f2_11}
M(\delta) = \left(
(E(\delta) x_i,x_j)_H
\right)_{i,j=0}^{N-1},\qquad \delta\in \mathfrak{B}( D_{\mathbf{r},\mathbf{l}} ),
\end{equation}
where $E$ is the spectral measure of $\mathcal{A}$.

\noindent
Then there exists a unitary transformation $V$, which maps $L^2(M)$ onto $H$ such that:
\begin{equation}
\label{f2_12}
V^{-1} S_j V = X_j,\qquad 1\leq j\leq \mathbf{r};
\end{equation}
\begin{equation}
\label{f2_13}
V^{-1} U_k V = W_k,\qquad 1\leq k\leq \mathbf{l}.
\end{equation}
Moreover, it holds:
\begin{equation}
\label{f2_14}
V \vec e_s = x_s,\qquad 0\leq s\leq N-1.
\end{equation}
\end{theorem}
In the case $\mathbf{r}=0$ relation~(\ref{f2_9}),(\ref{f2_12}) should be removed, and
in~(\ref{f2_10}) operators $S_j$ disappear. In the case $\mathbf{l}=0$ relation~(\ref{f2_13})
should be deleted and in~(\ref{f2_10}) operators $U_k$ disappear.
We are now ready to state and prove the following theorem.

\begin{theorem}
\label{t2_1}
Let $U_1, ...,U_n$ be pairwise commuting unitary operators in a separable Hilbert space $H$ ($n\in\mathbb{N}$). Then there exist
conjugations $J_1, ..., J_n$, $C$ in $H$, such that
\begin{equation}
\label{f2_30}
U_k = J_k C,\qquad k=1,...,n.
\end{equation}
Moreover, there exist conjugations $L_1, ..., L_n$, $K$ in $H$, such that
\begin{equation}
\label{f2_35}
U_k = K L_k,\qquad k=1,...,n.
\end{equation}
\end{theorem}
\textbf{Proof.} 
On the first step we shall represent the Hilbert space $H$ as a direct sum of some Hilbert spaces $H_j$. Each $H_j$ will reduce the
operators $U_k$, and $(U_1,...,U_n)$ will be a a commuting $SU$-set of order $(0,n)$ in $H_j$,
having the spectrum of multiplicity $1$. This can be done in a classical manner, like it was done for a single 
self-adjoint operator, see~\cite[Chapter 1, pp. 166-168]{cit_3500_Gelfand_Vilenkin_book}.
Namely, choose an arbitrary orthonormal basis $\mathfrak{F} = \{ f_j \}_{j=0}^\infty$ in $H$. Consider the following subspace:
$$ H_0 := \mathop{\rm span}\nolimits \{ U_1^{k_1} ... U_n^{k_n} f_0,\ k_1,...,k_n\in\mathbb{Z} \}. $$

On step $m$ ($m\in\mathbb{N}$) we choose the first element $u$ (if such exists, otherwise we stop the process)
of the basis $\mathfrak{F}$ which does not lie in $G_m := \oplus_{l=0}^{m-1} H_l$. 
We set $u' := u - P^H_{G_m} u$, and 
$$ H_m := \mathop{\rm span}\nolimits \{ U_1^{k_1} ... U_n^{k_n} u',\ k_1,...,k_n\in\mathbb{Z} \}. $$
We shall obtain the required decomposition:
\begin{equation}
\label{f2_39}
H = \oplus_{l=0}^{r} H_l,\qquad r\leq+\infty.
\end{equation}
Fix an arbitrary $l$: $0\leq l\leq r$.
Applying Theorem~\ref{t2_1_MN} for $(U_1,...,U_n)$ in the Hilbert space $H_l$, we conclude that
there exists a unitary transformation $V_l$, which maps $L^2(M_l)$ (with some scalar measure $M_l$ on 
$\mathfrak{B}( D_{0,n} )$) onto $H_l$ such that:
\begin{equation}
\label{f2_41}
V_l^{-1} U_{k;l} V_l = W_{k;l},\qquad 1\leq k\leq n,
\end{equation}
where
$$ W_{k;l} f(\varphi) = e^{i\varphi_k} f(\varphi),\qquad f\in L^2(M_l);\ 1\leq k\leq n, $$
and by $U_{k;l}$ we denoted the operator $U_k$ restricted to $H_l$.
Let us represent the operator $W_{k;l}$ as a product of two conjugations:
\begin{equation}
\label{f2_46}
W_{k;l} = J_{k;l} C_l,\qquad 1\leq k\leq n,
\end{equation}
where
\begin{equation}
\label{f2_48}
J_{k;l} f(\varphi) = e^{i\varphi_k} \overline{ f(\varphi) },\qquad 
C_l f(\varphi) = \overline{ f(\varphi) },\qquad
f\in L^2(M_l).  
\end{equation}
Define the following conjugations in $H_l$: 
$$ \widetilde J_{k;l} = V_l J_{k;l} V_l^{-1},\quad 1\leq k\leq n,\qquad \widetilde C_l = V_l C_l V_l^{-1}. $$
By~(\ref{f2_41}),(\ref{f2_48}) we conclude that
\begin{equation}
\label{f2_49}
U_{k;l} = \widetilde J_{k;l} \widetilde C_{l},\qquad 1\leq k\leq n.
\end{equation}
For an arbitrary $h = \sum_{l=0}^r h_l$, $h_l\in H_l$, $\sum_0^r \| h_l \|^2 < \infty$, we define the
following conjugations:
$$ \widetilde J_k h = \sum_{l=0}^r \widetilde J_{k;l} h_l,\quad 1\leq k\leq n;\quad
\widetilde C h = \sum_{l=0}^r \widetilde C_{l} h_l. $$
Then
$$ U_k = \widetilde J_k \widetilde C,\qquad 1\leq k\leq n, $$
and representation~(\ref{f2_30}) is verified. 
Representation~(\ref{f2_35}) can be obtained by applying~(\ref{f2_30})
for $(U_1^*,...,U_n^*)$.
$\Box$

A linear isometric operator $V$ in a Hilbert space $H$ is said to be \textit{essentially unitary}, if $\overline{V}$
is a unitary operator in $H$.
\begin{theorem}
\label{t2_3_m}
Let $\mathcal{T} = (A_1,...,A_\rho, B_1,...,B_\tau)$ ($\rho\in\mathbb{N}, \tau\in\mathbb{Z}_+, \rho + \tau\geq 2$) be 
a commuting tuple of symmetric and isometric operators (with a joint invariant dense domain $\mathcal{D}$) in 
a separable Hilbert space $H$.
Suppose that the following conditions hold:

\begin{itemize}
\item[(a)] The operators $A_2,...,A_\rho$ are essentially selfadjoint. The operators $B_1,...,B_\tau$
are essentially unitary. The closures of all of the above operators pairwise commute;

\item[(b)] $B_k \mathcal{D} = \mathcal{D}$, $k=1,...,\tau$;

\item[(c)] For some $z_0\in \mathbb{C}\backslash \mathbb{R}$, the operator $A_j$ ($2\leq j\leq\rho$)
restricted to the domain $(A_1-z_0 E_H)\mathcal{D}$ is essentially self-adjoint in a
Hilbert space
$(\overline{A_1} - z_0 E_H) D(\overline{A_1})$.

\item[(d)] There exists a conjugation $J$ in $H$ such that
$J \mathcal{D}\subseteq \mathcal{D}$, and
$$
A_j J = J A_j,\qquad j=1,...,\rho; $$
\begin{equation}
\label{f2_49_2}
B_k J = J B_k^{-1},\qquad k=1,...,\tau.
\end{equation}
\end{itemize}
Then there exists a self-adjoint operator $\widehat{A_1}\supseteq A_1$ in $H$, which commutes with $\overline{A_j}$ ($j=2,...,\rho$)
and with $\overline{B_k}$ ($k=1,...,\tau$).
\end{theorem}
\begin{remark}
\label{r2_1}
In the case $\rho = 1$ (or $\tau = 0$) the meaningless statements about operators $A_j$, $j\geq2$, in~(a),(c) (respectively about
operators $B_k$ in~(a),(b),(\ref{f2_49_2})) should be ignored in the conditions of the last theorem.
In the same way, the commutativity with non-existing operators should be ignored in the conclusion.
\end{remark}

\noindent
\textbf{Proof.} 
We shall use the constructions in a proof of Theorem~\ref{t1_1} (see~\cite{cit_9000_Zagorodnyuk_MFAT}) with necessary modifications and additions.
Consider the Cayley transformation of the operator $\overline{A_1}$:
$$
V_1 := (\overline{A_1} - \overline{z_0} E_H)(\overline{A_1}- z_0 E_H)^{-1}
= E_H + (z_0-\overline{z_0}) (\overline{A_1}- z_0 E_H)^{-1}, $$
and denote
\begin{equation}
\label{f2_49_5}
H_1 = (\overline{A_1} - z_0 E_H) D(\overline{A_1}),\
H_2 = H\ominus H_1,\
H_3 = (\overline{A_1} - \overline{z_0} E_H) D(\overline{A_1}),\
H_4 = H\ominus H_3. 
\end{equation}

If $\rho\geq 2$, then conditions of Theorem~\ref{t1_1} are satisfied with
$\mathbf{A}=A_1$, $\mathbf{B}=A_j$, $\mathbf{H}=H$, $\mathbf{J}=J$, for each $j\in\mathbb{Z}_{2,\rho}$.
Let us recall the construction in~\cite{cit_9000_Zagorodnyuk_MFAT} of a selfadjoint operator $\widetilde{A}_1 \supseteq A_1$, which
commutes with $\overline{A_j}$ in $H$ (for a fixed $j\in\mathbb{Z}_{2,\rho}$ at this moment).
Consider the Cayley transformation of the operator $\overline{A_j}$:
$$
U_j := (\overline{A_j} + i E_H)(\overline{A_j} - i E_H)^{-1}
= E_H + 2i (\overline{A_j} - iE_H)^{-1}.  $$
Denote by $U_{j,2}$ the restriction of $U_j$ to the subspace $H_2$, which is a unitary operator in $H_2$.
By the Godi\v{c}-Lucenko Theorem it can be decomposed as a product of two conjugations:
\begin{equation}
\label{f2_49_7}
U_{j,2} = K_j L_j.
\end{equation}
Then $\widetilde{A}_1$ is the inverse Cayley transformation of $V_1\oplus J K_j$.
If $K_j$ could be chosen independent of $j$, then $\widetilde{A}_1$ would commute with all $A_j$.
In the case $\tau = 0$, we can apply Theorem~\ref{t2_1} with $U_{j,2}$ ($2\leq j\leq\rho$)
which ensures that such a choice is indeed possible. Thus, for the case $\rho\geq 2,\tau=0$ the proof is complete.
In other cases we shall choose operators $K_j$ later, in order to take into account the operators $B_k$.

Suppose now that $\tau\geq 1$, and $\rho\in\mathbb{N}$.
For an arbitrary element $g = (A_1 - z_0 E_H) f$, $f\in\mathcal{D}$, we may write:
$$ B_k V_1 g = B_k (A_1 - \overline{z_0} E_H) f = (A_1 - \overline{z_0} E_H) B_k f; $$
$$ V_1 B_k g = V_1 (A_1 - z_0 E_H) B_k f = (A_1 - \overline{z_0} E_H) B_k f,\qquad k\in\mathbb{Z}_{1,\tau}. $$
By the continuity we get
\begin{equation}
\label{f2_50}
\overline{B_k} V_1 g = V_1 \overline{B_k} g,\qquad g\in H_1,\ k\in\mathbb{Z}_{1,\tau}. 
\end{equation}
By condition~(b) and the commutativity of $A_1$ and $B_k$ we obtain that
$$ B_k^{-1} A_1 = A_1 B_k^{-1},\qquad k\in\mathbb{Z}_{1,\tau}. $$
For an arbitrary element $g = (A_1 - z_0 E_H) f$, $f\in\mathcal{D}$, we may write:
$$ B_k^{-1} g = B_k^{-1} (A_1 - z_0 E_H) f = (A_1 - z_0 E_H) B_k^{-1} f\in H_1. $$
Therefore
$$ (\overline{B_k})^{-1} H_1\subseteq H_1. $$
In the same manner, we can see that
$$ \overline{B_k} H_1\subseteq H_1,\quad (\overline{B_k})^{-1} H_3\subseteq H_3,\quad \overline{B_k} H_3\subseteq H_3. $$
Therefore
$$ \overline{B_k} H_1 = H_1,\quad \overline{B_k} H_3 = H_3. $$
Thus, the operators $\overline{B_k}$ ($k\in\mathbb{Z}_{1,\tau}$) are unitary when restricted to each of the Hilbert spaces $H_i$, $i=1,2,3,4$.
Denote by $B_{k,2}$ the operator $\overline{B_k}$, restricted to $H_2$, $k\in\mathbb{Z}_{1,\tau}$.

In the case $\rho =1$, we may apply Theorem~\ref{t2_3_m} to the operators $B_{k,2}$, $k=1,...,\tau$, and obtain that
\begin{equation}
\label{f2_50_5}
B_{k,2} = \widetilde K R_k,\qquad k=1,...,\tau,
\end{equation}
with some conjugations $\widetilde K,R_k$ in $H_2$.
In the case $\rho\geq 2$, we apply Theorem~\ref{t2_3_m} to the operators $U_{j,2}$, $j\in\mathbb{Z}_{2,\rho}$;  $B_{k,2}$, $k\in\mathbb{Z}_{1,\tau}$,
and obtain that
\begin{equation}
\label{f2_50_7}
U_{j,2} = \widehat K G_j,\quad j\in\mathbb{Z}_{2,\rho},\quad B_{k,2} = \widehat K S_k,\quad k\in\mathbb{Z}_{1,\tau},
\end{equation}
with some conjugations $G_j, \widehat K,S_k$ in $H_2$.
We shall use the conjugation 
$$ K := 
\left\{ \begin{array}{cc} \widetilde K, & \mbox{if }\rho=1 \\
         \widehat K, & \mbox{if }\rho\geq 2 \end{array} \right., $$ 
in order to construct the desired extension of $A_1$.
As in the proof of Theorem~\ref{t1_1} (see formula~(4) in~\cite{cit_9000_Zagorodnyuk_MFAT}) we can check that
\begin{equation}
\label{f2_50_9}
\overline{A_1}^* J x = J \overline{A_1}^* x,\qquad x\in D(\overline{A_1}^*),
\end{equation}
and therefore
$$ J H_2 = H_4. $$
Set
\begin{equation}
\label{f2_57}
U_{2,4} = JK.
\end{equation}
The operator $U_{2,4}$ is linear and isometric, it maps $H_2$ onto $H_4$.
We choose $\widehat{A}_1$ to be the inverse Cayley's transform of $V_1\oplus U_{2,4}$.
Fix an arbitrary $k\in\mathbb{Z}_{1,\tau}$.
Let us check that $\widehat{A}_1$ commutes with the operator $\overline{B_k}$.
Taking into account relation~(\ref{f2_50}) we see, that it is enough to verify to the following condition:
\begin{equation}
\label{MFAT2f2_15_10}
\overline{B_k} U_{2,4} x = U_{2,4} \overline{B_k} x,\qquad x\in H_2.
\end{equation}
By the continuity and condition~(d) of the theorem we conclude that
\begin{equation}
\label{f2_59}
\overline{B_k} J = J (\overline{B_k})^{-1},\qquad k\in\mathbb{Z}_{1,\tau}. 
\end{equation}
Then
$$ U_{2,4} B_{k,2} U_{2,4}^{-1} y = J (B_{k,2})^{-1} J y = $$
$$ = J \overline{B_k}^{-1} J y = \overline{B_k} y,\qquad y\in H_4, $$
and relation~(\ref{MFAT2f2_15_10}) follows.
If $\rho=1$ this completes the proof. If $\rho\geq 2$, the considerations after~(\ref{f2_49_7}) show
that $\widehat A_1$ commutes with $A_j$ ($j\in\mathbb{Z}_{2,\rho}$) as well.
$\Box$

\begin{corollary}
\label{c2_1}
Let $\mathcal{T} = (A_1,...,A_\rho, B_1,...,B_\tau)$ ($\rho\in\mathbb{N}, \tau\in\mathbb{Z}_+$, $\rho + \tau\geq 2$) be 
a commuting tuple of symmetric and isometric operators (with a joint invariant dense domain $\mathcal{D}$) in 
a separable Hilbert space $H$.
Suppose that conditions~(b) and (d) of Theorem~\ref{t2_3_m} hold, as well as the following condition holds:
\begin{itemize}
\item[(a')] The operators $A_2,...,A_\rho$ are bounded. The operators $B_1,...,B_\tau$
are essentially unitary. 
\end{itemize}

Then there exists a self-adjoint operator $\widehat{A_1}\supseteq A_1$ in $H$, which commutes with $\overline{A_j}$ ($j=2,...,\rho$)
and with $\overline{B_k}$ ($k=1,...,\tau$).
\end{corollary}
\textbf{Proof.} 
It is clear that by the continuity condition~(a') implies condition~(a) of Theorem~\ref{t2_3_m}. 
Moreover, in our case condition~(c) of Theorem~\ref{t2_3_m} is trivially satisfied.
Thus, it remains
to apply Theorem~\ref{t2_3_m}.
$\Box$

We shall now illustrate our extension result by a concrete example.

\noindent
\textbf{The multidimensional power-trigonometric moment problem.} 
Fix arbitrary $r,l\in\mathbb{N}$.
For elements $(x,\varphi)\in D_{r,l}$ and $\mathbf{m} = (m_1,...,m_r)\in\mathbb{Z}_+^r$, $\mathbf{n} = (n_1,...,n_l)\in\mathbb{Z}^l$
($r,l\in\mathbb{N}$) we shall use the following multi-index notation:
$$ x^\mathbf{m} e^{i \mathbf{n}\varphi} = x_1^{m_1}...x_r^{m_r} e^{i n_1 \varphi_1}...e^{i n_l\varphi_l}. $$

We consider a moment problem which consists of
finding a non-negative measure $\mu$ on $\mathfrak{B}(D_{r,l})$ such that
\begin{equation}
\label{f2_63}
\int_{D_{r,l}} x^\mathbf{m} e^{i \mathbf{n}\varphi} d\mu = s_{\mathbf{m},\mathbf{n}},\qquad
\mathbf{m}\in\mathbb{Z}_+^r,\ \mathbf{n}\in\mathbb{Z}^l,
\end{equation}
where $\mathcal{S} := \{ s_{\mathbf{m},\mathbf{n}} \}_{\mathbf{m}\in \mathbb{Z}_+^r,\ \mathbf{n}\in \mathbb{Z}^l}$
is a prescribed set of complex numbers (moments).

\noindent
The moment problem~(\ref{f2_63}) is said to be \textit{the multidimensional Devinatz moment problem} or
\textit{the multidimensional power-trigonometric moment problem}. For the case $r=l=1$ see~\cite{cit_8950_Zagorodnyuk_JOT}
and references therein.

Suppose that the moment problem~(\ref{f2_63}) has a solution $\mu$.
Consider the following polynomial:
\begin{equation}
\label{f2_64}
p(x,\varphi) = \sum_{\mathbf{m}\in\mathbb{Z}_+^r,\ \mathbf{n}\in\mathbb{Z}^l} \alpha_{\mathbf{m},\mathbf{n}} x^{\mathbf{m}}
e^{i\mathbf{n}\varphi},\qquad \alpha_{\mathbf{m},\mathbf{n}}\in\mathbb{C},\ (x,\varphi)\in D_{r,l}, 
\end{equation}
where all but finite number of $\alpha_{\mathbf{m},\mathbf{n}}$s are zeros (such sequences of $\alpha_{\mathbf{m},\mathbf{n}}$s we shall
call \textit{finite}). 
Then
$$ 0\leq \int |p|^2 d\mu = \sum_{\mathbf{m},\mathbf{k}\in\mathbb{Z}_+^r,\ \mathbf{n},\mathbf{l}\in\mathbb{Z}^l}
\alpha_{\mathbf{m},\mathbf{n}} \overline{ \alpha_{\mathbf{k},\mathbf{l}} } \int x^{\mathbf{m}+\mathbf{k}} e^{i (\mathbf{n}-\mathbf{l} )\varphi} d\mu
= $$
$$ = \sum_{\mathbf{m},\mathbf{k}\in\mathbb{Z}_+^r,\ \mathbf{n},\mathbf{l}\in\mathbb{Z}^l}
\alpha_{\mathbf{m},\mathbf{n}} \overline{ \alpha_{\mathbf{k},\mathbf{l}} } s_{\mathbf{m}+\mathbf{k},\mathbf{n}-\mathbf{l}}.
$$
Therefore the following condition is necessary for the solvability of the moment problem~(\ref{f2_63}):
\begin{equation}
\label{f2_65}
\sum_{\mathbf{m},\mathbf{k}\in\mathbb{Z}_+^r,\ \mathbf{n},\mathbf{l}\in\mathbb{Z}^l}
\alpha_{\mathbf{m},\mathbf{n}} \overline{ \alpha_{\mathbf{k},\mathbf{l}} } s_{\mathbf{m}+\mathbf{k},\mathbf{n}-\mathbf{l}} \geq 0,
\end{equation}
for all finite sequences $\alpha_{\mathbf{m},\mathbf{n}}\in\mathbb{C}$.

Suppose now that the moment problem~(\ref{f2_63}) is given and condition~(\ref{f2_65}) holds.
Denote by $\mathfrak{L}$ a set of all polynomials of the form~(\ref{f2_64}). Notice that $\mathfrak{L}$
is a linear vector space with usual operations of addition and multiplication by a complex scalar (one needs to add the corresponding coefficients
by $x^{\mathbf{m}} e^{i\mathbf{n}\varphi}$
or multiply $\alpha_{\mathbf{m},\mathbf{n}}$s by a scalar).
Let $p$ be given by~(\ref{f2_64}) and
\begin{equation}
\label{f2_67}
q(x,\varphi) = \sum_{\mathbf{k}\in\mathbb{Z}_+^r,\ \mathbf{l}\in\mathbb{Z}^l} \beta_{\mathbf{k},\mathbf{l}} x^{\mathbf{k}}
e^{i\mathbf{l}\varphi},\qquad \beta_{\mathbf{k},\mathbf{l}}\in\mathbb{C},\ (x,\varphi)\in D_{r,l}, 
\end{equation}
where $\beta_{\mathbf{k},\mathbf{l}}$ is a finite sequence.
Consider the following functional:
\begin{equation}
\label{f2_70}
B(p,q) = 
\sum_{\mathbf{m},\mathbf{k}\in\mathbb{Z}_+^r,\ \mathbf{n},\mathbf{l}\in\mathbb{Z}^l}
\alpha_{\mathbf{m},\mathbf{n}} \overline{ \beta_{\mathbf{k},\mathbf{l}} } s_{\mathbf{m}+\mathbf{k},\mathbf{n}-\mathbf{l}}.
\end{equation}
The functional $B(p,q)$ is sesquilinear. By~(\ref{f2_65}) we see that $B(p,p)\geq 0$.  Taking the imaginary parts of 
$B(p+q,p+q)$ and $B(p+iq,p+iq)$, we conclude that $\overline{B(p,q)} = B(q,p)$.
We say that two polynomials $p,q\in\mathfrak{L}$ belong to the same equivalence class, if
$B(p-q,p-q)=0$. Making the completion in a set of the equivalence classes $\{ [p], p\in\mathfrak{L} \}$ we obtain a Hilbert space $H$.
It is clear that $H$ is separable (since one can consider polynomials with coefficients having rational real and imaginary parts).
Notice that
$$ ([p], [q])_H = B(p,q),\qquad p,q\in\mathfrak{L}. $$
Denote
\begin{equation}
\label{f2_75}
y_{\mathbf{m},\mathbf{n}} = [x^{\mathbf{m}} e^{i\mathbf{n}\varphi}],\qquad \mathbf{m}\in\mathbb{Z}_+^r,\ \mathbf{n}\in\mathbb{Z}^l.
\end{equation}
The following condition holds:
\begin{equation}
\label{f2_77}
(y_{\mathbf{m},\mathbf{n}}, y_{\mathbf{k},\mathbf{l}})_H = s_{\mathbf{m}+\mathbf{k},\mathbf{n}-\mathbf{l}},\qquad 
\mathbf{m,k}\in\mathbb{Z}_+^r,\ \mathbf{n,l}\in\mathbb{Z}^l.
\end{equation}
Denote $\mathcal{D} = \mathop{\rm Lin}\nolimits \{ y_{\mathbf{m},\mathbf{n}} \}_{\mathbf{m}\in\mathbb{Z}_+^r,\ \mathbf{n}\in\mathbb{Z}^l}$.
Observe that $\overline{\mathcal{D}}=H$. 
Set
$$ \vec e_{j;t} = (\delta_{1,j},\delta_{2,j},...,\delta_{t,j}),\qquad 1\leq j\leq t;\qquad t\in\mathbb{N}. $$
Define the following operators:
\begin{equation}
\label{f2_79}
A_j x = \sum_{\mathbf{m}\in\mathbb{Z}_+^r,\ \mathbf{n}\in\mathbb{Z}^l} \alpha_{\mathbf{m},\mathbf{n}} y_{\mathbf{m} +
\vec e_{j;r},\mathbf{n}},\qquad j\in\mathbb{Z}_{1,r}; 
\end{equation}
\begin{equation}
\label{f2_81}
B_k x = \sum_{\mathbf{m}\in\mathbb{Z}_+^r,\ \mathbf{n}\in\mathbb{Z}^l} \alpha_{\mathbf{m},\mathbf{n}} y_{\mathbf{m},
\mathbf{n}+ \vec e_{k;l}},\qquad k\in\mathbb{Z}_{1,l}; 
\end{equation}
where
\begin{equation}
\label{f2_83}
x = \sum_{\mathbf{m}\in\mathbb{Z}_+^r,\ \mathbf{n}\in\mathbb{Z}^l} \alpha_{\mathbf{m},\mathbf{n}} y_{\mathbf{m},\mathbf{n}}\in \mathcal{D}, 
\end{equation}
and $\alpha_{\mathbf{m},\mathbf{n}}$ is a finite sequence of complex numbers.
Let us check that these operators are well-defined.
Suppose that for an element $x\in\mathcal{D}$ with the representation~(\ref{f2_83}) there exists another representation:
\begin{equation}
\label{f2_85}
x = \sum_{\mathbf{m}\in\mathbb{Z}_+^r,\ \mathbf{n}\in\mathbb{Z}^l} \beta_{\mathbf{m},\mathbf{n}} y_{\mathbf{m},\mathbf{n}}, 
\end{equation}
where $\beta_{\mathbf{m},\mathbf{n}}$ is a finite sequence of complex numbers.
We need to show that
\begin{equation}
\label{f2_87}
\sum_{\mathbf{m}\in\mathbb{Z}_+^r,\ \mathbf{n}\in\mathbb{Z}^l} \alpha_{\mathbf{m},\mathbf{n}} y_{\mathbf{m} +
\vec e_{j;r},\mathbf{n}} =
\sum_{\mathbf{m}\in\mathbb{Z}_+^r,\ \mathbf{n}\in\mathbb{Z}^l} \beta_{\mathbf{m},\mathbf{n}} y_{\mathbf{m} +
\vec e_{j;r},\mathbf{n}},\qquad j\in\mathbb{Z}_{1,r}; 
\end{equation}
and
\begin{equation}
\label{f2_89}
\sum_{\mathbf{m}\in\mathbb{Z}_+^r,\ \mathbf{n}\in\mathbb{Z}^l} \alpha_{\mathbf{m},\mathbf{n}} y_{\mathbf{m},
\mathbf{n}+ \vec e_{k;l}} =
\sum_{\mathbf{m}\in\mathbb{Z}_+^r,\ \mathbf{n}\in\mathbb{Z}^l} \beta_{\mathbf{m},\mathbf{n}} y_{\mathbf{m},
\mathbf{n}+ \vec e_{k;l}},\qquad k\in\mathbb{Z}_{1,l}.
\end{equation}
Observe that for arbitrary $\mathbf{k}\in\mathbb{Z}_+^r, \mathbf{l}\in\mathbb{Z}^l$ we may write:
$$ \left(
\sum_{\mathbf{m}\in\mathbb{Z}_+^r,\ \mathbf{n}\in\mathbb{Z}^l} \alpha_{\mathbf{m},\mathbf{n}} y_{\mathbf{m} +
\vec e_{j;r},\mathbf{n}}, y_{\mathbf{k},\mathbf{l}}
\right)_H = 
\sum_{\mathbf{m}\in\mathbb{Z}_+^r,\ \mathbf{n}\in\mathbb{Z}^l} \alpha_{\mathbf{m},\mathbf{n}}
s_{\mathbf{m} + \vec e_{j;r} +\mathbf{k} ,\mathbf{n}-\mathbf{l}} = $$
$$ = \sum_{\mathbf{m}\in\mathbb{Z}_+^r,\ \mathbf{n}\in\mathbb{Z}^l} \alpha_{\mathbf{m},\mathbf{n}}
( y_{ \mathbf{m},\mathbf{n} }, y_{ \mathbf{k}+ \vec e_{j;r},\mathbf{l} } )_H 
= (x, y_{ \mathbf{k}+ \vec e_{j;r},\mathbf{l} })_H. $$
The same result we shall obtain if we take $\beta_{\mathbf{m},\mathbf{n}}$s instead of $\alpha_{\mathbf{m},\mathbf{n}}$s,
and relation~(\ref{f2_87}) follows.
On the other hand, we may write:
$$ \left\|
\sum_{\mathbf{m}\in\mathbb{Z}_+^r,\ \mathbf{n}\in\mathbb{Z}^l} (\alpha_{\mathbf{m},\mathbf{n}} - \beta_{\mathbf{m},\mathbf{n}})
y_{\mathbf{m}, \mathbf{n}+ \vec e_{k;l}}
\right\|^2 = $$
$$ = \sum_{\mathbf{m,k}\in\mathbb{Z}_+^r,\ \mathbf{n,l}\in\mathbb{Z}^l}
(\alpha_{\mathbf{m},\mathbf{n}} - \beta_{\mathbf{m},\mathbf{n}})
\overline{(\alpha_{\mathbf{k},\mathbf{l}} - \beta_{\mathbf{k},\mathbf{l}})}
(y_{\mathbf{m}, \mathbf{n}+ \vec e_{k;l}}, y_{\mathbf{k}, \mathbf{l}+ \vec e_{k;l}}) = $$
$$ = \sum_{\mathbf{m,k}\in\mathbb{Z}_+^r,\ \mathbf{n,l}\in\mathbb{Z}^l}
(\alpha_{\mathbf{m},\mathbf{n}} - \beta_{\mathbf{m},\mathbf{n}})
\overline{(\alpha_{\mathbf{k},\mathbf{l}} - \beta_{\mathbf{k},\mathbf{l}})}
(y_{\mathbf{m}, \mathbf{n}}, y_{\mathbf{k}, \mathbf{l}}) = $$
$$ = (x-x,x-x) = 0. $$
Thus, relation~(\ref{f2_89}) holds as well.
It is easily checked that all $A_j$ are symmetric, all $B_k$ are linear isometric, and all the above operators pairwise commute on
the joint invariant dense domain $\mathcal{D}$.
We shall need the following operator:
\begin{equation}
\label{f2_91}
J_0 x = \sum_{\mathbf{m}\in\mathbb{Z}_+^r,\ \mathbf{n}\in\mathbb{Z}^l} \overline{\alpha_{\mathbf{m},\mathbf{n}}} y_{\mathbf{m},
-\mathbf{n}},
\end{equation}
where $x$ is given by relation~(\ref{f2_83}).
It is easily checked that $J_0$ is well-defined and
$$ (J_0 x, J_0 y) = (y,x),\qquad x,y\in\mathcal{D}. $$
By the continuity we extend $J_0$ to a conjugation $J$ on $H$.

If $r\geq 2$ we additionally suppose that the following condition holds:

\noindent
(B): \textit{
For some $j_0\in\mathbb{Z}_{1,r}$,
there exist positive constants $C_j$ ($j\in\mathbb{Z}_{1,r}\backslash\{ j_0 \}$) such that
$$ \sum_{\mathbf{m,k}\in\mathbb{Z}_+^r,\ \mathbf{n,l}\in\mathbb{Z}^l}
\alpha_{\mathbf{m},\mathbf{n}} \overline{ \alpha_{\mathbf{k},\mathbf{l}} }
s_{ \mathbf{m} + \mathbf{k} + 2\vec e_{j;r}, \mathbf{n} - \mathbf{l} } \leq $$
\begin{equation}
\label{f2_93}
\leq C_j
\sum_{\mathbf{m,k}\in\mathbb{Z}_+^r,\ \mathbf{n,l}\in\mathbb{Z}^l}
\alpha_{\mathbf{m},\mathbf{n}} \overline{ \alpha_{\mathbf{k},\mathbf{l}} }
s_{ \mathbf{m} + \mathbf{k}, \mathbf{n} - \mathbf{l} },
\end{equation}
for all finite sequences of complex numbers $\alpha_{\mathbf{m},\mathbf{n}}$ and all $j\in\mathbb{Z}_{1,r}\backslash\{ j_0 \}$.
}

Condition~(B) ensures that operators $A_j$ ($j\in\mathbb{Z}_{1,r}\backslash\{ j_0 \}$) are bounded. 
By Corollary~\ref{c2_1} we conclude that there exists a self-adjoint extension $\widetilde A_{j_0}\supseteq A_{j_0}$
in $H$, which commutes with the closures of the rest of $A_j$s and $B_k$s. 

By the induction one can check that
$$ y_{(m_1,...,m_r),(n_1,...,n_l)} = A_1^{m_1} ... A_r^{m_r} B_1^{n_1} ... B_l^{n_l} y_{\mathbf{0}_r,\mathbf{0}_l}, $$
\begin{equation}
\label{f2_95}
\mathbf{m}=(m_1,...,m_r)\in\mathbb{Z}_+^r,\ \mathbf{n}=(n_1,...,n_l)\in\mathbb{Z}^l,
\end{equation}
where $\mathbf{0}_t$ means the null vector of size $t$.
Then
$$ s_{\mathbf{m},\mathbf{n}} = ( y_{(m_1,...,m_r),(n_1,...,n_l)}, y_{\mathbf{0}_r,\mathbf{0}_l})_H =
( A_1^{m_1} ... A_r^{m_r} B_1^{n_1} ... B_l^{n_l} y_{\mathbf{0}_r,\mathbf{0}_l}, y_{\mathbf{0}_r,\mathbf{0}_l})_H = $$
$$ = \int_{D_{r,l}} x^\mathbf{m} e^{i \mathbf{n}\varphi} d\widetilde\mu,\quad
\mathbf{m}=(m_1,...,m_r)\in\mathbb{Z}_+^r,\ \mathbf{n}=(n_1,...,n_l)\in\mathbb{Z}^l, $$
where $\widetilde\mu(\delta) = (E(\delta) y_{\mathbf{0}_r,\mathbf{0}_l}, y_{\mathbf{0}_r,\mathbf{0}_l})_H$, $\delta\in\mathfrak{B}(D_{r,l})$,
and $E(\delta)$ is the spectral measure for the commuting SU-set of
$\overline{A_j}$ ($j\in\mathbb{Z}_{1,r}\backslash\{ j_0 \}$), $\widetilde A_{j_0}$ and $\overline{B_k}$ ($k\in\mathbb{Z}_{1,l}$).
Thus, $\widetilde\mu$ is a solution of the moment problem~(\ref{f2_63}).

\begin{theorem}
\label{t2_4}
Let the moment problem~(\ref{f2_63}) be given with some moments 
$\{ s_{\mathbf{m},\mathbf{n}} \}_{\mathbf{m}\in \mathbb{Z}_+^r,\ \mathbf{n}\in \mathbb{Z}^l}$.
If $r=1$, then condition~(\ref{f2_65}) is necessary and sufficient for the solvability of the moment problem~(\ref{f2_63}).

If $r\geq 2$, then condition~(\ref{f2_65}) is still necessary, while
condition~(\ref{f2_65}) together with condition~(B) are sufficient
for the solvability of the moment problem~(\ref{f2_63}).
\end{theorem}
\textbf{Proof.} 
All statements of the theorem were already proved in the preceding considerations.
$\Box$

\begin{center}
{\large\bf Extensions of commuting tuples of symmetric and isometric operators.}
\end{center}
\begin{center}
{\bf S.M. Zagorodnyuk}
\end{center}

In this paper we study extensions of commuting tuples of symmetric and isometric operators to
commuting tuples of self-adjoint and unitary operators. Some conditions which ensure the existence of
such extensions are presented. A multidimensional analog of the Godi\v{c}-Lucenko Theorem is proved.
An application to a multidimensional power-trigonometric moment problem is given.

\vspace{1cm}

V. N. Karazin Kharkiv National University \newline\indent
School of Mathematics and Computer Sciences \newline\indent
Department of Higher Mathematics and Informatics \newline\indent
Svobody Square 4, 61022, Kharkiv, Ukraine

Sergey.M.Zagorodnyuk@gmail.com; Sergey.M.Zagorodnyuk@univer.kharkov.ua

}
\end{document}